# INTEGRAL TRANSFORM OF THE GALUĖ TYPE STRUVE FUNCTION

D.L. Suthar[1], S.D. Purohit[2] and K.S. Nisar[3]

[1]Department of Mathematics, Wollo University, Dessie Campus
P.O. Box: 1145, South Wollo, Amhara Region, Ethiopia
Email: dlsuthar@gmail.com

[2]Department of HEAS (Mathematics)
Rajasthan Technical University, Kota, India
Email: sunil_a_purohit@yahoo.com

[3]Department of Mathematics, College of Arts and Science-Wadi Al-dawasir
Prince Sattam bin Abdulaziz University, Saudi Arabia.
Email: ksnisar1@gmail.com

This paper refers to the study of generalized Struve type function. Using generalized Galuė type Struve function (GTSF) by Nisar et al. [13], we derive various integral transform, including Euler transform, Laplace transform, Whittakar transform, K-transform and Fractional Fourier transform. The results are expressed in terms of generalized Wright function. The transform established here are general in nature and are likely to find useful in applied problem of Sciences, engineering and technology. Certain interesting special cases of the main results are also considered.

## 1. INTRODUCTION

Integral transforms have been widely used in various problems of mathematical physics and applied mathematics. This information has inspired the study of several integral transforms with varity of special functions. The present paper deals with the evaluation of the Euler transform, Laplace transform, Whittakar transform, K-transform and Fractional Fourier transform of the Galuė type Struve function recently introduced by [13]. This work is motivated to the authors for their recent work [12, 13, 14].

**Definition 1.1 Generalized Galuė type Struve function**

Recently, Nisar et. al. [13], defined as following generalized form of Struve function named as generalized Galuė type Struve function (GTSF) as:

$${}_a W_{p,b,c,\xi}^{\lambda,\mu}(z) = \sum_{k=0}^{\infty} \frac{(-c)^k}{\Gamma(\lambda k + \mu)\Gamma\left(ak + \frac{p}{\xi} + \frac{b+2}{2}\right)} \left(\frac{z}{2}\right)^{2k+p+1}, \quad a \in \mathbb{N}\ p,b,c \in \mathbb{C} \quad (1.1)$$

where $\lambda > 0$, $\xi > 0$ and $\mu$ is an arbitrary parameter.

For the definition of the struve function and its more generalization, the interested reader may refer to the many papers ( Bhow-mick [3, 4], Kanth [6], Singh [17, 18], Nisar and Atangana [11], Singh [19].

**Special cases**

When $\lambda = a = 1$, $\mu = 3/2$ and $\xi = 1$ in equation (1.1), it reduces to generalization of Struve function which is defined by Orhan and Yagmur [15, 16 ].

$$H_{p,b,c}(z) = \sum_{k=0}^{\infty} \frac{(-c)^k}{\Gamma\left(k+\frac{3}{2}\right)\Gamma\left(k+p+\frac{b+2}{2}\right)} \left(\frac{z}{2}\right)^{2k+p+1}, \quad p,b,c \in \mathbb{C} \tag{1.2}$$

Details related to the function $H_{p,b,c}(z)$ and its particular cases can be seen in Barics [1, 2], Mondal and Swaminathan [10], Mondal and Nisar [9].

Throught this paper, we need the following well known facts and rules.

**Definition 1.2 Euler Transform** (Sneddon [20])

The Euler transform of a function $f(z)$ is defined as

$$B\{f(z); a,b\} = \int_0^1 z^{a-1}(1-z)^{b-1} f(z)dz \qquad a,b \in \mathbb{C}, \ \Re(a) > 0, \ \Re(b) > 0 \tag{1.3}$$

**Definition 1.3 Laplace Transform** (Sneddon [20])

The Laplace transform of a function $f(t)$, denoted by $F(s)$, is defined by the equation

$$F(s) = (Lf)(s) = L\{f(t); s\} = \int_0^{\infty} e^{-st} f(t)dt \qquad \Re(s) > 0 \tag{1.4}$$

Provided the integral (1.4) is convergent and that the function $f(t)$, is continous for $t>0$ and of exponential order as $t \to \infty$, (1.4) may be symbolically written as

$$F(s) = L\{f(t); s\} \quad \text{or} \quad f(t) = L^{-1}\{F(s); t\} \tag{1.5}$$

**Definition 1.4 Whittakar Transform** (Whittakar and Watson [22])

$$\int_0^{\infty} e^{-\frac{1}{2}t} t^{\zeta-1} W_{\tau,\omega}(t)dt = \frac{\Gamma\left(\frac{1}{2}+w+\zeta\right)\Gamma\left(\frac{1}{2}-w+\zeta\right)}{\Gamma(1-\tau+\zeta)} \tag{1.6}$$

where $\Re(w \pm \zeta) > -1/2$ and $W_{\tau,\omega}(t)$ is the Whittakar confluent hypergeometric function

$$W_{\omega,\zeta}(z) = \frac{\Gamma(-2\omega)}{\Gamma\left(\frac{1}{2}-\tau-\omega\right)} M_{\tau,\omega}(z) + \frac{\Gamma(2\omega)}{\Gamma\left(\frac{1}{2}+\tau+\omega\right)} M_{\tau,-\omega}(z) \tag{1.7}$$

where $M_{\tau,\omega}(z)$ is defined by

$$M_{\tau,\omega}(z) = z^{1/2+\omega} e^{-1/2 z} {}_1F_1\left(\frac{1}{2}+\omega-\tau; \ 2\omega+1; \ z\right) \tag{1.8}$$

**Definition 1.5 K-Transform** (Erdělyi et al.[5])

This transform is defined by the following integral equation

$$\Re_v[f(x); p] = g[p; v] = \int_0^{\infty} (px)^{1/2} K_v(px) f(x) dx \tag{1.9}$$

where $\Re(p) > 0$; $K_v(x)$ is the Bessel function of the second kind defined by [5, p.332]

$$K_v(z) = \left(\frac{\pi}{2z}\right)^{1/2} W_{0,v}(2z)$$

where $W_{0,\upsilon}(.)$ is the Whittakar function defined in equation (1.7).

The following result given in Mathai et al. [8, p. 54, eq. 2.37] will be used in evaluating the integrals:

$$\int_0^\infty t^{\rho-1} K_\upsilon(ax)dx = 2^{\rho-2} a^{-\rho} \Gamma\left(\frac{\rho\pm\upsilon}{2}\right); \quad \Re(a)>0; \Re(\rho\pm\upsilon)>0. \tag{1.10}$$

**Definition 1.6 Fractional Fourier Transform** (Luchko et al.[7])

The fractional Fourier transform of order $\alpha$, $0<\alpha\leq 1$ is defined by

$$\hat{u}_\alpha(\omega) = \mathfrak{J}_\alpha[u](\omega) = \int_R e^{i\omega^{(1/\alpha)}t} u(t)dt \tag{1.11}$$

when $\alpha=1$, equation (1.11) reduces to the conventional Fourier transforms and for $\omega>0$, it reduces to the Fractional Fourier transform defined by Luchko et al. [7].

The compositions are expressed in terms of the generaliged Wright hypergeometric function ${}_p\psi_q(z)$ (see, for detail, Srivastava and Karson [21]), for $z\in\mathbb{C}$ complex, $a_i, b_j \in\mathbb{C}$ and $\alpha_i,\beta_j \in\Re$, where $(\alpha_i,\beta_j \neq 0;$ $i=1,2,...,p; j=1, 2,...,q)$, is defined as below:

$$_p\psi_q(z) = {}_p\psi_q\left[\begin{matrix}(a_i,\alpha_i)_{1,p}\\(b_j,\beta_j)_{1,q}\end{matrix}\bigg|z\right] = \sum_{k=0}^\infty \frac{\prod_{i=1}^p \Gamma(a_i+\alpha_i k) z^k}{\prod_{j=1}^q \Gamma(b_j+\beta_j k) k!}. \tag{1.12}$$

Introduced by Wright [23], the generalized Wright function (1.12) and proved several theorems on the asymptotic expansion of ${}_p\psi_q(z)$ (for instance, see [23, 24, 25]) for all values of the argument $z$, under the condition:

$$\sum_{j=1}^q \beta_j - \sum_{i=1}^p \alpha_i > -1.$$

## 2. INTEGRAL TRANSFORMS OF ${}_a w_{p,b,c,\xi}^{\lambda,\mu}(z)$

**Theorem 2.1 (Euler Transform)** Let $a\in\mathbb{N}$ $p,b,c,r,s \in\mathbb{C}$; be such that

$$B\left\{{}_a w_{p,b,c,\xi}^{\lambda,\mu}\left(x^{1/2} z\right); r, s\right\} = \left(\frac{\sqrt{x}}{2}\right)^{p+1} \Gamma(s) \, {}_2\psi_3\left[\begin{matrix}(p+r+1, 2), (1, 1);\\(\mu, \lambda), \left(\frac{p}{\xi}+\frac{b}{2}+1, a\right), (p+r+s+1, 2);\end{matrix}\bigg|\frac{-cx}{4}\right] \tag{2.1}$$

where $\Re(r)>0, \Re(s)>0$, $\lambda>0$, $\xi>0$ and $\mu$ is an arbitrary parameter.

**Proof.** Using (1.1) and (1.3), it gives

$$B\left\{{}_a w_{p,b,c,\xi}^{\lambda,\mu}\left(x^{1/2} z\right); r, s\right\} = \int_0^1 z^{r-1}(1-z)^{s-1} {}_a w_{p,b,c,\xi}^{\lambda,\mu}\left(x^{1/2} z\right) dz$$

$$= \sum_{k=0}^\infty \frac{(-c)^k}{\Gamma(\lambda k+\mu)\Gamma\left(ak+\frac{p}{\xi}+\frac{b+2}{2}\right)} \left(\frac{x^{1/2}}{2}\right)^{2k+p+1} \int_0^1 z^{2k+p+r+1-1}(1-z)^{s-1} dz$$

$$= \sum_{k=0}^{\infty} \frac{(-c)^k}{\Gamma(\lambda k + \mu)\Gamma\left(ak + \frac{p}{\xi} + \frac{b+2}{2}\right)} \left(\frac{x^{1/2}}{2}\right)^{2k+p+1} B(p+r+1+2k,\ s)$$

$$= \left(\frac{x^{1/2}}{2}\right)^{p+1} \sum_{k=0}^{\infty} \frac{1}{\Gamma(\lambda k + \mu)\Gamma\left(\frac{p}{\xi} + \frac{b+2}{2} + ak\right)} \frac{\Gamma(p+r+1+2k)\,\Gamma(s)}{\Gamma(p+r+s+1+2k)} \left(\frac{(-cx)^{1/2}}{2}\right)^k$$

In accordance with the definition of (1.12), we obatain the result (2.1). This completes the proof of the theorem.

**Corollary 2.1** For $\lambda = a = 1$, $\mu = 3/2$ and $\xi = 1$, equation (2.1) reduces in the following form

$$B\left\{H_{p,b,c}\left(x^{1/2}z\right);\ r,\ s\right\} = \left(\frac{x^{1/2}}{2}\right)^{p+1} \Gamma(s)\ _2\Psi_3\left[\begin{array}{c}(p+r+1,\ 2),\ (1,\ 1);\\ \left(p+\frac{b}{2}+1,\ 1\right),\ (p+r+s+1,\ 2),\ (3/2,\ 1);\end{array}\left|\frac{-cx}{4}\right.\right] \qquad (2.2)$$

**Theorem 2.2 (Laplace Transform)** Let $a \in \mathbb{N}$ $p,b,c \in \mathbb{C}$; be such that

$$L\left\{_a w_{p,b,c,\xi}^{\lambda,\mu}\left(x^{1/2}z\right);\ s\right\} = \left(\frac{\sqrt{x}}{2}\right)^{p+1} s^{-(p+2)}\ _2\Psi_3\left[\begin{array}{c}(p+2,\ 2),\ (1,\ 1);\\ (\mu,\ \lambda), \left(\frac{p}{\xi}+\frac{b}{2}+1,\ a\right);\end{array}\left|\frac{-cx}{4s^k}\right.\right] \qquad (2.3)$$

where $\lambda > 0$, $\xi > 0$ and $\mu$ is an arbitrary parameter.

**Proof.** Using (1.1) and (1.4), it gives

$$L\left\{_a w_{p,b,c,\xi}^{\lambda,\mu}\left(x^{1/2}z\right);\ s\right\} = \int_0^{\infty} e^{-sz}\ _a w_{p,b,c,\xi}^{\lambda,\mu}\left(x^{1/2}z\right)dz$$

$$= \sum_{k=0}^{\infty} \frac{(-c)^k}{\Gamma(\lambda k + \mu)\Gamma\left(ak + \frac{p}{\xi} + \frac{b+2}{2}\right)} \left(\frac{x^{1/2}}{2}\right)^{2k+p+1} \int_0^{\infty} z^{2k+p+1-1} e^{zs}\ dz$$

$$= \left(\frac{x^{1/2}}{2}\right)^{p+1} \sum_{k=0}^{\infty} \frac{\Gamma(p+2+2k)}{\Gamma(\lambda k + \mu)\Gamma\left(ak + \frac{p}{\xi} + \frac{b+2}{2}\right)} \left(\frac{-cx^{1/2}}{2}\right)^{2k} L\left\{\frac{z^{p+2+2k-1}}{\Gamma(p+2+2k)};\ s\right\}$$

In accordance with the definition of (1.12), we obatain the result (2.3). This completes the proof of the theorem.

**Corollary 2.2** For $\lambda = a = 1$, $\mu = 3/2$ and $\xi = 1$, equation (2.3) reduces in the following form

$$L\left\{H_{p,b,c}\left(x^{1/2}z\right);\ s\right\} = \left(\frac{\sqrt{x}}{2}\right)^{p+1} s^{-(p+2)}\ _2\Psi_3\left[\begin{array}{c}(p+2,\ 2),\ (1,\ 1);\\ \left(p+\frac{b}{2}+1,\ 1\right),\ (3/2,\ 1);\end{array}\left|\frac{-cx}{4s^k}\right.\right] \qquad (2.4)$$

**Theorem 2.3 (Whittakar Transform)** Let $a \in \mathbb{N}$ $p,b,c \in \mathbb{C}$; $\Re(\zeta) > 0$, $\Re(w \pm \zeta) > -1/2$ be such that

$$\int_0^{\infty} t^{\zeta-1} e^{-\frac{t}{2}} W_{\tau,\omega}(t)\ _a w_{p,b,c,\xi}^{\lambda,\mu}\left(x^{1/2}t\right)dt = \left(\frac{\sqrt{x}}{2}\right)^{p+1}\ _3\Psi_3\left[\begin{array}{c}(w+\zeta+p+3/2,\ 2),\ (-w+\zeta+p+3/2,\ 2),\ (1,\ 1);\\ (\mu,\ \lambda), \left(\frac{p}{\xi}+\frac{b}{2}+1,\ a\right),\ (-\tau+\zeta+p+2,\ 2);\end{array}\left|\frac{-cx}{4}\right.\right] \qquad (2.5)$$

where $\Re(e) > |\Re(\omega)| - 1/2$, $\Re(p) > 0$, $\lambda > 0$, $\xi > 0$ and $\mu$ is an arbitrary parameter.

**Proof.** Using (1.1) and (1.6), it gives

$$\int_0^\infty t^{\zeta-1} e^{-\frac{t}{2}} W_{\tau,\omega}(t)_a w_{p,b,c,\xi}^{\lambda,\mu}\left(x^{1/2} t\right) dt = \int_0^\infty t^{\zeta-1} e^{-\frac{t}{2}} W_{\tau,\omega}(t) \sum_{k=0}^\infty \frac{(-c)^k}{\Gamma(\lambda k + \mu)\Gamma\left(ak + \frac{p}{\xi} + \frac{b+2}{2}\right)} \left(\frac{x^{1/2} t}{2}\right)^{2k+p+1} dt$$

$$= \sum_{k=0}^\infty \frac{(-c)^k}{\Gamma(\lambda k + \mu)\Gamma\left(ak + \frac{p}{\xi} + \frac{b+2}{2}\right)} \left(\frac{x^{1/2}}{2}\right)^{2k+p+1} \int_0^\infty e^{-\frac{t}{2}} t^{\zeta+p+2k+1-1} W_{\tau,\omega}(t) dt$$

$$= \sum_{k=0}^\infty \frac{(-c)^k}{\Gamma(\lambda k + \mu)\Gamma\left(ak + \frac{p}{\xi} + \frac{b+2}{2}\right)} \left(\frac{x^{1/2}}{2}\right)^{2k+p+1} \frac{\Gamma\left(\frac{1}{2} + w + \zeta + p + 2k + 1\right)\Gamma\left(\frac{1}{2} - w + \zeta + p + 2k + 1\right)}{\Gamma(1 - \tau + \zeta + p + 2k + 1)}$$

$$= \left(\frac{x^{1/2}}{2}\right)^{p+1} \sum_{k=0}^\infty \frac{\Gamma\left(\frac{1}{2} + w + \zeta + p + 2k + 1\right)\Gamma\left(\frac{1}{2} - w + \zeta + p + 2k + 1\right)}{\Gamma(\lambda k + \mu)\,\Gamma(1 - \tau + \zeta + p + 2k + 1)\,\Gamma\left(ak + \frac{p}{\xi} + \frac{b+2}{2}\right)} \left(\frac{(-cx)^{1/2}}{2}\right)^{2k}$$

In accordance with the definition of (1.12), we obatain the result (2.5). This completes the proof of the theorem.

**Corollary 2.3** For $\lambda = a = 1$, $\mu = 3/2$ and $\xi = 1$, equation (2.5) reduces in the following form

$$\int_0^\infty t^{\zeta-1} e^{-\frac{t}{2}} W_{\tau,\omega}(t) H_{p,b,c}\left(x^{1/2} t\right) dt = \left(\frac{\sqrt{x}}{2}\right)^{p+1} {}_3\Psi_3\left[\begin{array}{c}(w+\zeta+p+3/2,\ 2),\ (-w+\zeta+p+3/2,\ 2),\ (1,\ 1);\\ \left(\frac{3}{2},\ 1\right), \left(p+\frac{b}{2}+1,\ 1\right),\ (-\tau+\zeta+p+2,\ 2);\end{array}\bigg|\frac{-cx}{4}\right] \quad (2.6)$$

**Theorem 2.4 (K-Transform)** Let $a \in \mathbb{N}$ $p,b,c,\rho \in \mathbb{C}$; be such that

$$\int_0^\infty t^{\rho-1} K_\upsilon(\omega t)_a w_{p,b,c,\xi}^{\lambda,\mu}\left(x^{1/2} t\right) dt = 2^{\rho+p-1} \omega^{1-\rho-p} \left(\frac{\sqrt{x}}{2}\right)^{p+1} {}_3\Psi_2\left[\begin{array}{c}\left(\frac{\rho+p+\upsilon+1}{2},1\right), \left(\frac{\rho+p-\upsilon+1}{2},1\right),\ (1,\ 1);\\ (\mu,\ \lambda), \left(\frac{p}{\xi}+\frac{b}{2}+1,\ a\right);\end{array}\bigg|\frac{-cx}{\omega^2}\right]$$

(2.7)

where $\Re(\omega) > 0$, $\Re(\rho \pm \upsilon) > 0$. $\lambda > 0$, $\xi > 0$ and $\mu$ is an arbitrary parameter.

**Proof.** Using (1.1) and (1.10), it gives

$$\int_0^\infty t^{\rho-1} K_\upsilon(\omega t) H_{p,b,c}\left(x^{1/2} t\right) dt = \int_0^\infty t^{\rho-1} K_\upsilon(\omega t) \sum_{k=0}^\infty \frac{(-c)^k}{\Gamma(\lambda k + \mu)\Gamma\left(ak + \frac{p}{\xi} + \frac{b+2}{2}\right)} \left(\frac{x^{1/2} t}{2}\right)^{2k+p+1} dt$$

$$= \sum_{k=0}^\infty \frac{(-c)^k}{\Gamma(\lambda k + \mu)\Gamma\left(ak + \frac{p}{\xi} + \frac{b+2}{2}\right)} \left(\frac{x^{1/2}}{2}\right)^{2k+p+1} \int_0^\infty t^{\rho+p+2k+1-1} K_\upsilon(\omega t) dt$$

$$= \sum_{k=0}^{\infty} \frac{(-c)^k}{\Gamma(\lambda k + \mu)\Gamma\left(ak + \frac{p}{\xi} + \frac{b+2}{2}\right)} \left(\frac{x^{1/2}}{2}\right)^{2k+p+1} 2^{\rho+p+2k-1} \omega^{1-\rho-p-2k} \Gamma\left(\frac{(\rho+p+2k+1)\pm\upsilon}{2}\right)$$

$$= 2^{\rho+p-1} \omega^{1-\rho-p} \left(\frac{\sqrt{x}}{2}\right)^{p+1} \sum_{k=0}^{\infty} \frac{\Gamma\left(\frac{\rho+p+\upsilon+1}{2}+k\right)\Gamma\left(\frac{\rho+p-\upsilon+1}{2}+k\right)}{\Gamma(\lambda k + \mu)\Gamma\left(ak + \frac{p}{\xi} + \frac{b+2}{2}\right)} \left(\frac{(-cx)^{1/2}}{\omega}\right)^{2k}$$

In accordance with the definition of (1.12), we obatain the result (2.7). This completes the proof of the theorem.

**Corollary 2.4** For $\lambda = a = 1$, $\mu = 3/2$ and $\xi = 1$, equation (2.7) reduces in the following form

$$\int_0^{\infty} t^{\rho-1} K_{\upsilon}(\omega t) H_{p,b,c}\left(x^{1/2} t\right) dt = 2^{\rho+p-1} \omega^{1-\rho-p} \left(\frac{\sqrt{x}}{2}\right)^{p+1} {}_3\Psi_2 \left[\begin{array}{c} \left(\frac{\rho+p+\upsilon+1}{2},1\right), \left(\frac{\rho+p-\upsilon+1}{2},1\right), (1,\,1); \\ \left(\frac{3}{2},1\right), \left(p+\frac{b}{2}+1,\,1\right); \end{array} \middle| \frac{-cx}{\omega^2} \right]$$

(2.8)

### 3. FRACTIONAL FOURIER TRANSFORMS (FFT) OF ${}_a w^{\lambda,\mu}_{p,b,c,\xi}(z)$

**Theorem 3.1 (Euler Transform)** Let $a \in \mathbb{N}$ $p,b,c \in \mathbb{C}$; be such that

$$\mathfrak{I}_{\zeta}\left[{}_a w^{\lambda,\mu}_{p,b,c,\xi}\left(x^{1/2} t\right)\right] = \left(\frac{\sqrt{x}}{2}\right)^{p+1} \sum_{k=0}^{\infty} \frac{1}{(i)^{2k+p+2} \omega^{(2k+p+2)/\zeta} (-1)^{2k+p+1}} \frac{\Gamma(2k+p+2)}{\Gamma(\lambda k + \mu)\Gamma\left(ak + \frac{p}{\xi} + \frac{b+2}{2}\right)} \left(\frac{(-cx)^{1/2}}{2}\right)^{2k}$$

(3.1)

where $\zeta > 0$, $\lambda > 0$, $\xi > 0$ and $\mu$ is an arbitrary parameter.

**Proof.** Using (1.1) and (1.11), it gives

$$\mathfrak{I}_{\zeta}\left[{}_a w^{\lambda,\mu}_{p,b,c,\xi}\left(x^{1/2} t\right)\right](\omega) = \int_R e^{i\omega^{(1/\zeta)} t} \sum_{k=0}^{\infty} \frac{(-c)^k}{\Gamma(\lambda k + \mu)\Gamma\left(ak + \frac{p}{\xi} + \frac{b+2}{2}\right)} \left(\frac{x^{1/2} t}{2}\right)^{2k+p+1} dt$$

$$= \sum_{k=0}^{\infty} \frac{(-c)^k}{\Gamma(\lambda k + \mu)\Gamma\left(ak + \frac{p}{\xi} + \frac{b+2}{2}\right)} \left(\frac{x^{1/2}}{2}\right)^{2k+p+1} \int_R e^{i\omega^{(1/\zeta)} t} t^{2k+p+1} dt$$

If we set $i\omega^{(1/\zeta)} t = \eta$, then

$$= \sum_{k=0}^{\infty} \frac{(-c)^k}{\Gamma(\lambda k + \mu)\Gamma\left(ak + \frac{p}{\xi} + \frac{b+2}{2}\right)} \left(\frac{x^{1/2}}{2}\right)^{2k+p+1} \int_{-\infty}^{0} e^{-\eta} \left(\frac{-\eta}{i\omega^{1/\zeta}}\right)^{2k+p+1} \left(\frac{-d\eta}{i\omega^{1/\zeta}}\right)$$

$$= \sum_{k=0}^{\infty} \frac{(-c)^k}{\Gamma(\lambda k + \mu)\Gamma\left(ak + \frac{p}{\xi} + \frac{b+2}{2}\right)} \left(\frac{x^{1/2}}{2}\right)^{2k+p+1} \frac{1}{(i)^{2k+p+2} \omega^{(2k+p+2)/\zeta} (-1)^{2k+p+1}} \int_0^{\infty} e^{-\eta} \eta^{2k+p+1} d\eta$$

$$= \sum_{k=0}^{\infty} \frac{(-c)^k}{\Gamma(\lambda k + \mu)\Gamma\left(ak + \frac{p}{\xi} + \frac{b+2}{2}\right)} \left(\frac{x^{1/2}}{2}\right)^{2k+p+1} \frac{\Gamma(2k+p+2)}{(i)^{2k+p+2} \omega^{(2k+p+2)/\zeta}(-1)^{2k+p+1}}$$

$$= \left(\frac{\sqrt{x}}{2}\right)^{p+1} \sum_{k=0}^{\infty} \frac{1}{(i)^{2k+p+2} \omega^{(2k+p+2)/\zeta}(-1)^{2k+p+1}} \frac{\Gamma(2k+p+2)}{\Gamma(\lambda k + \mu)\Gamma\left(ak + \frac{p}{\xi} + \frac{b+2}{2}\right)} \left(\frac{(-cx)^{1/2}}{2}\right)^{2k}$$

This completes the proof of the theorem.

**Corollary 3.1** For $\lambda = a = 1$, $\mu = 3/2$ and $\xi = 1$, equation (3.1) reduces in the following form

$$\mathfrak{I}_\zeta \left[ H_{p,b,c}\left(x^{1/2} t\right) \right] = \left(\frac{\sqrt{x}}{2}\right)^{p+1} \sum_{k=0}^{\infty} \frac{1}{(i)^{2k+p+2} \omega^{(2k+p+2)/\zeta}(-1)^{2k+p+1}} \frac{\Gamma(2k+p+2)}{\Gamma\left(k+\frac{3}{2}\right)\Gamma\left(k+P+\frac{b+2}{2}\right)} \left(\frac{(-cx)^{1/2}}{2}\right)^{2k}$$